\documentclass[a4paper,12pt,reqno]{amsart}
\usepackage{latexsym}
\usepackage{csquotes}
\usepackage{amsmath,amsthm,amssymb}%,amscd
\usepackage{fullpage}
\usepackage{xcolor}
\usepackage{parskip}
\usepackage{mathtools}
\usepackage{thmtools}
\mathtoolsset{showonlyrefs}

\usepackage[pdfencoding=unicode,pdftex,bookmarks=true,bookmarksnumbered]{hyperref}
\definecolor{citecolour}{rgb}{0.0, 0.0, 0.8}
\definecolor{urlcolour}{rgb}{1,0.5,0}
\colorlet{linkcolour}{green!50!black}
\hypersetup{colorlinks,breaklinks,
		linkcolor=linkcolour,citecolor=citecolour,
		filecolor=linkcolour, urlcolor=urlcolour}

\theoremstyle{plain}
\newtheorem{prevtheorem}{Theorem}

\newtheorem{theorem}{Theorem}[section]
\newtheorem{proposition}[theorem]{Proposition}
\newtheorem{lemma}[theorem]{Lemma}
\newtheorem{corollary}[theorem]{Corollary}

\theoremstyle{definition}

\numberwithin{equation}{section}
\usepackage[shortlabels]{enumitem}
\begingroup
 \makeatletter
 \@for\theoremstyle:=definition,remark,plain\do{%
 \expandafter\g@addto@macro\csname th@\theoremstyle\endcsname{%
 \addtolength\thm@preskip\parskip
 }%
 }

\DeclareMathOperator{\Hom}{Hom}

\DeclareMathOperator{\ccr}{\mathcal{T}}
\newcommand{\gensub}[1]{\left\langle{#1}\right\rangle}

\newcommand{\R}{\mathbb{\R}}

\renewcommand{\leq}{\leqslant}
\renewcommand{\geq}{\geqslant}

\newenvironment{proofof}{{\bf {Proof.} }}{\hfill $\blacksquare$ \\}
\newenvironment{proofofA}{{\bf {Proof of Theorem~\ref{Thm:A}.}}}{\hfill $\blacksquare$ \\}
\newenvironment{proofofB}{{\bf {Proof of Theorem~\ref{Thm:B}.}}}{\hfill $\blacksquare$ \\}
\renewcommand{\phi}{\varphi}

\begin{document}

\title{Common transversals and complements\\ in abelian groups}
\date{\today}
\author{S. Aivazidis}
\address{Department of Mathematics and Applied Mathematics, University of Crete, Voutes Campus, 70013 Heraklion, Greece}
\email{s.aivazidis@uoc.gr}

\author{M. Loukaki}
\address{Department of Mathematics and Applied Mathematics, University of Crete, Voutes Campus, 70013 Heraklion, Greece}
\email{mloukaki@uoc.gr}

\author{B. Sambale}
\address{Institut f\"{u}r Algebra, Zahlentheorie und Diskrete Mathematik, Leibniz Universit\"{a}t Hannover, Welfengarten 1, 30167 Hannover, Germany}
\email{sambale@math.uni-hannover.de}

\thanks{The first author is supported by the Hellenic Foundation for Research and Innovation, Project HFRI-FM17-1733. The third author is supported by the German Research Foundation (\mbox{SA 2864/4-1}).}

\renewcommand{\sectionautorefname}{Section}

\begin{abstract}
Given a finite abelian group $G$ and cyclic subgroups $A$, $B$, $C$ of $G$
of the same order, we find necessary and sufficient conditions for $A$, $B$, $C$
to admit a common transversal for the cosets they afford.
For an arbitrary number of cyclic subgroups we give a sufficient criterion when there exists a common complement. 
Moreover, in several cases where a common transversal exists, we provide concrete constructions.
\end{abstract}



\maketitle
\tableofcontents

\section{Introduction}
Suppose that $G$ is a group 
with subgroups $A$, $B$ of the same finite index. 
Then Hall's celebrated Marriage Theorem guarantees 
that we can find a \emph{common} transversal for both $A$, $B$. 
That is, there exists a set $T$ whose elements comprise 
a complete and non-redundant set of coset representatives for both $A$ and $B$. Furthermore,  the total number of such common transversals has been computed in~\cite{ALM-transversal}. 
If three or more subgroups are involved, common transversals may not exist. For instance, the three subgroups of order $2$ in the Klein group do not possess such a common transversal. 

The recent interest in seeking a common transversal 
for several subgroups of an abelian group 
stems from a question of Steinhaus~\cite{sierpinski1958probleme} in the 1950s. 
Steinhaus asked if there exists a subset of the plane
which can \emph{tile} the plane when \emph{translated} 
by any one of the lattices that arise 
by rotating the integer lattice around the origin. 
In the language used in this paper, 
Steinhaus asked if there exists a common transversal of all lattices $R_\theta\mathbb{Z}^2$, 
where $\theta \in [0, 2\pi)$ and $R_\theta$ denotes rotation by $\theta$ around the origin.

It was only proved in this century~\cite{jackson2002lattice} 
that the answer is indeed affirmative. 
The variant of the problem where the subset of the plane 
is asked to be Lebesgue measurable 
(but tiling is only demanded almost everywhere) 
has sparked much more interest and is still open. 
The best results to date for the measurable problem can be found in~\cite{kolountzakis1999steinhaus}. 

Variations of the Steinhaus problem have taken many forms, 
but the one most relevant to this paper was first studied in~\cite{kolountzakis1997multi}; 
where the question was posed if we can find a subset of the plane 
which is a common transversal for a finite number of lattices in the plane. 
Surprisingly this problem has an affirmative solution in the measurable sense 
when the duals of the finite set of lattices has a direct product 
(the dual of the lattice $A\mathbb{Z}^d$ is the lattice $A^{-\top}\mathbb{Z}^d$). 
In~\cite{kolountzakis1997multi} the problem was first posed of 
when a finite set of subgroups of an abelian group of the same finite index 
admit a common transversal and it was proved~\cite[Thm.~1]{kolountzakis1997multi} 
that if the subgroups $A_1, \ldots, A_n$ of $G$ are direct factors in $G$, 
then they always admit a common transversal in $G$.

In this paper we are interested in tackling a specific case
of the general problem outlined above. 
In particular, we will have something useful to say when cyclic
subgroups of a finite abelian group are involved and a common
transversal is sought. 
As mentioned earlier, such a transversal does not always exist.
Recall that if $K \leq G$, with $G$ being arbitrary,
then $K$ is said to have a \emph{complement} in $G$
(or to be complemented in $G$) in case there exists a subgroup $H$
such that $G = KH$ and $K \cap H = 1$.
Observe that such a complement $H$, if it exists, is 
a transversal of $K$ in $G$ that moreover inherits the group structure of the parent group.

The first of our two main theorems reads as follows.

\begin{prevtheorem}\label{Thm:A} 
Let $A_1,\ldots,A_t$ be complemented isomorphic subgroups of a finite abelian group $G$. If the smallest prime divisor of $|A_1|$ is at least $t$, then $A_1,\ldots,A_t$ have a common complement in $G$.
\end{prevtheorem}

\autoref{Thm:A} implies that two isomorphic complemented subgroups in abelian groups always have a common complement. This is false for non-abelian groups as can be seen in the dihedral group of order $8$.
The proof of \autoref{Thm:A} accompanied with more detailed statements will be given in \autoref{Sec:Complements}.

Our second main theorem provides a complete description
of the situation as regards transversals when three cyclic subgroups are involved.

\begin{prevtheorem}\label{Thm:B} 
Let $G$ be a finite abelian group with cyclic subgroups $A$, $B$, $C$ of the same order. 
Then $A$, $B$, $C$ do not share a common transversal in $G$ 
if and only if $A$ (and thus $B$ and $C$) has even order 
and the product $A_2B_2C_2$ of their Sylow $2$-subgroups satisfies
\[
A_2B_2C_2/I = A_2/I \times B_2 /I= A_2/I \times C_2 /I = B_2/I \times C_2/I\,,
\]
where $I \coloneqq A_2\cap B_2\cap C_2$.
\end{prevtheorem}

When it comes to more than three subgroups, larger primes play a role. For instance the elementary abelian group of order $p^2$ is the union of $p+1$ subgroups of order $p$, but there cannot be a common transversal. 

\section{Reductions}\label{Sec:Reductions}

We outline below some notational
conventions that we will use throughout the paper.
\begin{itemize}
\item $[n]\coloneqq\{1,2,\ldots,n\}$.
\item A cyclic group of order $n$ is denoted by $C_n$, while $S_n$ is the symmetric group of degree $n$.
\item $G$ will always denote a finite abelian group.
\item A group is called homoyclic if it is the direct product of isomorphic cyclic groups.
\item For a prime $p$, $G_p$ denotes the unique Sylow $p$-subgroup of $G$.
\item If $H \leq G$ and $a, b \in G$ we write $a \equiv b \pmod{H}$ if $a H = b H$. 
\item For an integer $n$ let $\Gamma_n(G) \coloneqq \langle g\in G:g^n=1\rangle$. If $G$ is a $p$-group, we use the standard notation $\Omega(G) \coloneqq \Gamma_p(G)$. 
\item For $A_1,\ldots,A_n\leq G$ let $\ccr_G(A_1,\ldots,A_n)$ be the set of common transversals of $A_1,\ldots,A_n$ in $G$. 
Similarly, let $\mathcal{X}_G(A_1,\ldots,A_n)$ be the set of common complements of $A_1,\ldots,A_n$ in $G$.
\end{itemize}

We first observe that some (easy) reductions can be made.

\begin{lemma}\label{lemma:SubgTrans}
Let $A\leq H\leq G$ be finite abelian groups. 
If $X\in\ccr_G(A)$, then $X\cap H\in\ccr_H(A)$. 
\end{lemma}
\begin{proofof}
For any coset $hA$ of $A$ in $H$ there exists $x \in X$ so that $xA = hA$, 
since $X$ is a transversal of $A$ in $G$ and $hA$ is a coset of $A$ in $G$ as well.
Therefore, $x \in H$ and the lemma follows.
\end{proofof}

\begin{lemma}\label{Prop:GisProd}
For subgroups $A_1,\ldots,A_n\leq G$ we have
$\ccr_G(A_1,\ldots,A_n) \neq \emptyset$ if and only if 
$\ccr_{A_1\ldots A_n}(A_1,\ldots,A_n) \neq \emptyset$.
\end{lemma}
\begin{proofof}
If $\ccr_G(A_1,\ldots,A_n)\ne\emptyset$, then $\ccr_{A_1\ldots A_n}(A_1,\ldots,A_n)\ne\emptyset$ by \autoref{lemma:SubgTrans}.

Conversely, assume that $S \in \ccr_{A_1\ldots A_n}(A_1,\ldots,A_n)$ 
with $S = \{s_1, s_2, \ldots, s_m\}$ 
and let $T = \{t_1, t_2, \ldots, t_k\}$
be a transversal for $A_1\ldots A_n$ in $G$. 
Clearly the set 
$ST = \{ s_i t_j \mid i \in [m] , \, j \in [k]\}$ 
is a transversal for $A_1,\ldots,A_n$ in $G$ 
and thus $\ccr_G(A_1,\ldots,A_n)\ne\emptyset$.
\end{proofof}

\begin{lemma}\label{Prop:IntersectionIsTriv}
Let $A_1,\ldots,A_n\leq G$ and $N\leq A_1\cap\ldots\cap A_n$. Then $\{g_1, \ldots, g_m\} \in \ccr_G(A_1,\ldots,A_n)$ if and only if $\{g_1N, \ldots, g_mN\} \in \ccr_{G/N}(A_1/N,\ldots, A_n/N)$.
In particular, $\ccr_G(A_1,\ldots,A_n)\ne\emptyset$ if and only if $\ccr_{G/N}(A_1/N,\ldots, A_n/N)\ne\emptyset$.
\end{lemma}

\begin{proofof}
This follows easily from $|G:A_i|=|G/N:A_i/N|$ and 
the equivalence between
$g_i\equiv g_j\pmod{A}$ and $g_iN\equiv g_jN\pmod{A/N}$.
\end{proofof}

\noindent 
The next result shows that to decide whether or not
a number of  subgroups of an abelian group possess a
common complement in the parent group,
it suffices to decide whether their Sylow $p$-subgroups
are complemented, or equivalently to assume
that $G$ has prime power order. 

\begin{lemma}\label{Lem:Bijection}
Let $G$ be an abelian group of order $|G| = p_1^{a_1} \cdots p_r^{a_r}$ 
and let $A$ be a subgroup of $G$ that has a complement in $G$. 
Then there is a canonical bijection $\mathcal{X}_G(A)\to \bigtimes_{i=1}^r \mathcal{X}_{G_{p_i}}(A_{p_i})$.
\end{lemma}

\begin{proofof}
To begin with, observe that each subgroup of $G$ (thus also $G$ itself)
is the direct product of its Sylow $p$-subgroups.
Since $G_p$ is unique for each prime $p$,
we have $A_p = A \cap G_p\leq G_p$.
Now let 
$f : \mathcal{X}_G(A) \to \bigtimes_{i=1}^r \mathcal{X}_{G_{p_i}}(A_{p_i})$ 
be given by the rule
\[
\mathcal{X}_G(A) \ni H \mapsto \left(H \cap G_{p_1}, \ldots, H \cap G_{p_r}\right).
\]
We will argue that $f$ is 1-1 and onto.
To see that it is onto, 
let $(H_1, \ldots, H_r) \in\bigtimes_{i=1}^r \mathcal{X}_{G_{p_i}}(A_{p_i})$
and put $H = \prod_{i=1}^r H_i$.
Then $H\in \mathcal{X}_G(A)$ by order considerations 
and $f(H) = (H_1, \ldots, H_r)$.
On the other hand, let $H, K \leq \mathcal{X}_G(A)$ and suppose that
$f(H) = f(K)$.
Then
\[
\left(H \cap G_{p_1}, \ldots, H \cap G_{p_r}\right) =
\left(K \cap G_{p_1}, \ldots, K \cap G_{p_r}\right).
\]
By our previous remark each Sylow subgroup of $H$
coincides with the corresponding Sylow subgroup of $K$.
Since both $H$ and $K$ are the internal direct
products of their Sylow subgroups,
it follows that $H = K$ establishing the
desired injectivity of $f$.
\end{proofof}

The preceding lemma clearly implies:
\begin{corollary}\label{Prop:ReduceToSylow_complements}
Let $G$ be abelian group and $A_1, \ldots, A_n\leq G$.
Then $A_1, \ldots, A_n$ share a common complement in $G$
if and only if their Sylow $p$-subgroups $A_{1_p}, \ldots, A_{n_p}$  share a common complement in $G_p$ for all prime divisors $p$ of $|G|$.
\end{corollary}

\section{Common complements}\label{Sec:Complements}

Recall that the fundamental theorem of finite abelian groups asserts that $G$ has a primary decomposition $G=G_1\times\ldots\times G_n$ where each $G_i$ is a cyclic group of prime-power order. 
The usual proof is based on the fact that cyclic subgroups of maximal order have complements. 
The following consequence of the Krull--Remak--Schmidt theorem gives an exact criterion.

\begin{lemma}\label{KRS}
Let $G = G_1 \times \ldots \times G_n$ be the primary decomposition of an abelian group $G$. 
Then $A \leq G$ has a complement in $G$ if and only if there exists a subset $I\subseteq[n]$ 
such that the projection $\pi_I : A \to \bigtimes_{i\in I} G_i$ is an isomorphism.
\end{lemma}

\begin{proofof}
Suppose first that $\pi_I$ is an isomorphism. 
Let $B \coloneqq \bigtimes_{i\in[n]\setminus I} G_i$. 
For $a \in A\cap B$ we have $\pi_I(a) = 1$ and therefore $a = 1$. 
Since $|AB|=\prod_{i\in I}|G_i|\prod_{i\in[n]\setminus I}|G_i|=|G|$, 
$B$ is a complement of $A$ in $G$.

Now assume conversely that $G=A\times B$. Let 
\[
G=A_1\times\ldots\times A_r\times B_1\times\ldots\times B_s
\]
be the primary decomposition. 
By the Krull--Remak--Schmidt theorem (see~\cite[Satz~I.12.3]{Huppert}), 
there exists $J \subseteq [n]$ such that $G = A_1 \times \ldots \times A_r \times \bigtimes_{j \in J}G_j$. 
Clearly, $\pi_I$ with $I \coloneqq [n]\setminus J$ is an isomorphism. 
\end{proofof}

This is already sufficient to prove \autoref{Thm:A}.

\medskip
\begin{proofofA}
By \autoref{Prop:ReduceToSylow_complements} we may assume that $G$ is a $p$-group where $p\geq t$. 
We argue by induction on the rank of $A_1\cong\ldots\cong A_t$. Suppose first that $A_i$ is not cyclic. Let $A_i=A_{i1}\times A_{i2}$ be a non-trivial decomposition such that $A_{1j}\cong\ldots\cong A_{tj}$ for $j=1,2$. By hypothesis, $A_i$ has a complement $K_i$ in $G$. It is easy to see that $A_{i2}\times K_i$ is a complement of $A_{i1}$ in $G$. Hence, by induction there exists a common complement $H$ of $A_{11},\ldots,A_{t1}$ in $G$. Let $B_i:=A_i\cap H$. Since $A_{i1}B_i=A_{i1}H\cap A_i=A_i$ and $A_{i1}\cap B_i\le A_{i1}\cap H=1$, $A_{i1}$ is a complement of $B_i$ in $A_i$. Consequently, $A_{i1}K_i$ is a complement of $B_i$ in $G$. Finally, $A_{i1}K_i\cap H$ is a complement of $B_i$ in $H$. Hence, by induction there exists a common complement $K$ of $B_1,\ldots,B_t$ in $H$. Now we have $A_iK=A_iB_iK=A_iH=G$ and $A_i\cap K=A_i\cap H\cap K=B_i\cap K=1$. So $K$ is also a complement of $A_1,\ldots,A_t$ in $G$.

For the remainder of the proof, we may assume that $A_1,\ldots,A_t$ are cyclic. Let $G=G_1\times\ldots\times G_n$ be the primary decomposition. By \autoref{KRS}, there exist $i_1,\ldots,i_t\in[n]$ such that the projection $\pi_{i_j}:A_j\to G_{i_j}$ is an isomorphism for $j=1,\ldots,t$. In particular, $|G_{i_1}|=\ldots=|G_{i_t}|$. Let $H\coloneqq G_{i_1}\ldots G_{i_t}\cong A_1^s$ for some $s\leq t$ (note that the $G_{i_j}$ are not necessarily distinct). Let $G=H\times K$ and let $\pi_H:G\to H$ be the projection to $H$. By construction, the restriction of $\pi_H$ to each $A_i$ is injective.
Let $B_i \coloneqq \pi_H(A_i)\cap\Omega(H)$ for $i\in[t]$. Since $|\Omega(H)|=p^s$, $\Omega(H)$ has exactly $\frac{p^s-1}{p-1}$ maximal subgroups. Each $B_i$ can only be contained in the preimages of the $\frac{p^{s-1}-1}{p-1}$ maximal subgroups of $\Omega(H)/B_i$.
Since 
\[\frac{p^s-1}{p^{s-1}-1}>p\geq t,\] 
there must be one maximal subgroup of $\Omega(H)$ not containing any $B_i$. Consequently, we find a subgroup $L\leq H$ such that $L\cong A_1^{s-1}$ and $L\cap\pi_H(A_i)=1$ for $i=1,\ldots,t$ (the case $s=1$ with $L=1$ is allowed). 
Suppose that $a\in A_i\cap (L\times K)$. Then $\pi_H(a)\in L\cap\pi_H(A_i)=1$ and $a=1$ since $\pi_H$ is injective on $A_i$. Moreover, 
\[
|A_i||L\times K|=|A_i|^s|K|=|HK|=|G|.
\]
Hence, $L\times K$ is a common complement of $A_1, \ldots, A_t$ in $G$. 
\end{proofofA}

The proof above shows that the following stronger statement holds for cyclic subgroups of $p$-groups.

\begin{corollary}\label{cor:pgrps}
Let $A_1,\ldots,A_t$ be complemented cyclic subgroups of the same order of a finite abelian $p$-group $G$. If there are at most $p$ distinct subgroups $\Omega(A_1),\ldots,\Omega(A_t)$, then $A_1,\ldots,A_t$ have a common complement in $G$.
\end{corollary}

Note that the bound $t \leq p$ in the preceding corollary is sharp
owing to the existence of the elementary abelian group of order $p^2$ 
and its $p+1$ distinct subgroups of order $p$.

Our next result addresses the existence of a common transversal.
\begin{corollary}\label{cor:homocyclic}
Let $A_1,\ldots,A_t$ be homocyclic subgroups of the same order 
of a finite abelian group $G$. 
If the smallest prime divisor of $|A_1|$ is at least $t$, 
then $A_1, \ldots, A_t$ have a common transversal in $G$.
\end{corollary}

\begin{proofof}
Suppose that $K$ is the subgroup generated by $A_1, \ldots, A_t$.
We will show that $A_1, \ldots, A_t$ admit a common complement in $K$.
First, we argue that each subgroup $A_i$ is complemented in $K$
and we observe that it will suffice to prove the claim for $A_1$.
Note that since $K = A_1 \ldots A_t$ and
the $A_i$'s are homocyclic, it follows that $\exp(K) = \exp(A_1)$. 
Now we argue by induction on the rank of $A_1$,
where the base case is true since then $A_1$ is cyclic
of maximal order. Let $C$ be cyclic of maximal order in $A_1$.
Then $C$ has a complement, say $H$, in $K$ and so $K = C \times H$.
By Dedekind's lemma we have $A_1 = C \times (A_1 \cap H)$,
where $A_1 \cap H$ is homocyclic of rank one less than the rank of $A_1$.
Also, the exponent of $A_1 \cap H$ is equal to the exponent of $H$
and so the induction hypothesis applies to $A_1 \cap H$ in $H$ and
ensures the existence of a complement for $A_1 \cap H$ in $H$, say $D$.
Then it is easy to see that $D$ is a complement for $A_1$ in $K$, as claimed.
Now~\autoref{Thm:A} applies to the collection $A_1, \ldots, A_t$ in $K$
proving that there is a common complement. Thus $A_1, \ldots, A_t$ have
a common transversal in $K$ and so by~\autoref{Prop:GisProd} they
admit a common transversal in $G$.
\end{proofof}

A direct consequence of the preceding corollary
is that three homocyclic subgroups of the same odd order
always admit a common transversal in any parent group.

To obtain more refined results, we now start to count common complements. 
For this purpose the following basic fact is useful.

\begin{lemma}\label{numcompl}
Suppose that the abelian group $G$ has a complemented subgroup $A$. 
Then the number of complements is 
\[
|\mathcal{X}_G(A)| =  \left\lvert \Hom(G/A,A) \right\rvert  =  \left\lvert \Hom(A,G/A) \right\rvert .
\]
\end{lemma}
\begin{proofof}
Let $T\in\mathcal{X}_G(A)$ be fixed. For every $U\in\mathcal{X}_G(A)$ and every $t\in T$ there exists a unique element $\tau_U(t)\in A$ such that $\tau_U(t)\equiv t \pmod{U}$. It is easy to see that $\tau_U : T\to A$ is a homomorphism and $U$ is uniquely determined by $\tau_U$. Conversely, every homomorphism $\tau : T\to A$ defines a complement of $A$ 
as $U \coloneqq \{\tau(t)t^{-1}:t\in T\}$. 
In particular, $|\mathcal{X}_G(A)| =  \left\lvert  \Hom(T,A)  \right\rvert = \left\lvert \Hom(G/A,A) \right\rvert $. The last equality is a general fact from duality.
\end{proofof}

Note that there is a natural isomorphism $\Hom(A\times B,C)\cong\Hom(A,C)\times\Hom(B,C)$. In some situations this cardinality simplifies considerably. 

\begin{corollary}\label{corNrComp}
Suppose that $A\cong C_{n_1}\times\ldots\times C_{n_k}$ has a complement in $G$. Then 
\[
|\mathcal{X}_G(A)| = |\Gamma_{n_1}(G/A)|\ldots |\Gamma_{n_k}(G/A)|.
\]
\end{corollary}
\begin{proofof}
Let $A=\langle a_1\rangle\times\ldots\times\langle a_k\rangle$ with $|\langle a_i\rangle|=n_i$ for $i=1,\ldots,k$. Then every homomorphism $f:A\to G/A$ is uniquely determined by $f(a_i)\in\Gamma_{n_i}(G/A)$ for $i=1,\ldots,k$ and conversely, each such choice leads in fact to a homomorphism. 
\end{proofof}

\begin{lemma}\label{Lem:ProdA_i}
Let $G=A_1\times \cdots \times  A_t\times B$ be an abelian group. Then 
there are exactly 
\[
\prod_{i=2}^t\lvert\mathrm{Iso}(A_1,A_i)\rvert  \left\lvert \Hom(A_1,B) \right\rvert 
\]
common complements for $A_1,\ldots,A_t$ in $G$, where $\mathrm{Iso}(A_1,A_i)$ is the set of isomorphisms $A_1\to A_i$. 
\end{lemma}

\begin{proofof}
If $X\in\mathcal{X}_G(A_1,\ldots,A_t)$, then $A_1\cong G/X\cong A_i$ for each $i\geq 1$. Hence, we may assume that $A_1\cong\ldots\cong A_t$.
 Obviously, 
 \[
 T \coloneqq \bigtimes_{i=2}^tA_i \times B
 \]
 is a complement of $A_1$ in $G$. 
 To every complement $U$ of $A_1$ in $G$ there exists the corresponding homomorphism $\tau_U:T\to A_1$ such that $U= \{\tau_U(t)t^{-1}:t\in T\}$, 
% %$\tau_U(t)\equiv t \pmod{U}$ for every $t\in T$
 as in \autoref{numcompl}. Now fix $2\leq i\leq t$ and note  that $U$ can only be complement of $A_i$ too, if the restriction of $\tau_U$ to $A_i$ is injective, the reason being that $\tau_U(t)t^{-1} \in A_i$ for some $t \in T$ if and only if $\tau_U(t)=1 $ and $t \in A_i$.  But then $\tau_U(A_i)=A_1$ since $A_1\cong A_i$.
 In this case, $\tau_U$ decomposes into a product of isomorphisms $A_i\to A_1$ and a homomorphism $B\to A_1$. Conversely, it can be checked that each such $\tau$ (that is  the product of $t-1$ isomomorphisms from $A_i $ to $A$ along with a homomorphism from $B$ to $A_1$)  defines a common complement $\{\tau(t)t^{-1}:t\in T\}$ of
 $A_1,\ldots,A_t$. 
\end{proofof}

For cyclic groups the above clearly implies:
\begin{corollary}\label{Cor:cor1}
Let $G=A_1\times \cdots \times  A_t\times B$ with cyclic subgroups $A_1,\ldots,A_t$ of order $n$. Then 
\[
|\mathcal{X}_G(A_1,\ldots,A_t)|=\phi(n)^{t-1}|\Gamma_n(B)|,
\]
where $\phi$ is the totient function.
\end{corollary}

Now we restrict our attention to cyclic subgroups of maximal order (which are always complemented).

\begin{proposition}\label{Prop:NumberCommonComplementAB}
Let $A$, $B$ be cyclic subgroups of maximal order of an abelian $p$-group $G$.
Then
\[|\mathcal{X}_G(A,B)|=\begin{cases}
|G:A|&\text{if }A\cap B>1,\\
\phi(|G:A|)&\text{if }A\cap B=1.
\end{cases}\]
\end{proposition}

\begin{proofof}
Assume first that $A \cap B > 1$. 
Then every complement of $A$ is also a complement of $B$ and vice versa. 
To see this, observe that if $T$ is a complement of $A$ but not of $B$ 
then $T \cap B \neq 1$ and thus $\Omega(B) \leq T\cap B \cap A$, 
which is clearly a contradiction. Hence 
$\mathcal{X}_G(A,B) = \mathcal{X}_G(A)$ and 
\autoref{corNrComp} implies that $\mathcal{X}_G(A) = |G : A|$. 
So the proposition holds in this case.

\noindent 
Assume now that $A \cap B = 1$ and thus $AB = A \times B= C_{p^n}^2$ where $|A|=p^n$. Since by hypothesis, $p^n$ is the maximal element order in $G$, it follows easily from \autoref{KRS} that $AB$ has a complement in $G$ say $G = A \times B \times S$. 
We apply \autoref{Cor:cor1} to get exactly 
$\phi(p^n) |\Gamma_{p^n}(S)|$ common complements of $A$ and $B$ in $G$. 
Observe that $\Gamma_{p^n}(S)= S$ and thus 
\[
\phi(p^n) |\Gamma_{p^n}(S)|= \phi(p^n) \frac{|G|}{p^{2n}}= \phi(|G|/p^n)= \phi(|G : A|).
\]
The proof is complete.
\end{proofof}

\begin{corollary}
Let $G$ be an abelian group 
and let $A$, $B$ be cyclic subgroups of $G$ of maximal order
and index $s$ in $G$. 
Then the proportion of complements 
of $A$ in $G$ that are simultaneously complements for $B$ in $G$ is at least
$\phi(s)/s$.
\end{corollary}

\begin{proofof}
The desired proportion is 
\[
\frac{|\mathcal{X}_G(A,B)|}{|\mathcal{X}_G(A)|}. 
\]
If we write $n_{p}$ for the number of
common complements of the Sylow $p$-subgroups $A_{p}$ and $B_{p}$ in $G_{p}$, 
then according to \autoref{Lem:Bijection} we have 
\[
|\mathcal{X}_G(A,B)|=\prod_{p\, \mid\, |G|} n_{p}.
\]
But $n_{p}$ equals $|G_{p} : A_{p}|$, 
if $A_{p} \cap B_{p} > 1$ 
and $\phi(|G_{p} : A_{p}|)$, 
if $A_{p} \cap B_{p} = 1$ 
by \autoref{Prop:NumberCommonComplementAB}.
Hence in all cases we have $n_{p} \geq \phi(|G_{p} : A_{p}|)$
and thus the number of common complements of $A$, $B$ in $G$
is at least 
\[
\prod_{p \, \mid\,  |G|} \phi(|G_{p} : A_{p}|) = 
\phi\Bigl(\prod_{p \, \mid\,  |G|} |G_{p} : A_{p}|\Bigr) = 
\phi(|G : A|) = \phi(s).
\]

Similarly for $A$ we get 
\[
|\mathcal{X}_G(A)|= |\Gamma_{|A|}(G/A)|=|G/A|=s
\]
complements in $G$ by \autoref{Cor:cor1}. Therefore, 
the desired proportion is at least $\phi(s)/s$, as wanted.
\end{proofof}

Notice that for each prime power index $s=p^a$ the proportion 
of common complements is at least $1 - 1/p \geq 1/2$,
but for general $n$ there is no positive lower bound
since 
\[
\liminf_{n \to \infty} \frac{\phi(n)}{n} = 0
\]
(consider $n$ a product of distinct primes).

We can now give a quantitative version of \autoref{Thm:A}.

\begin{theorem}\label{Thm:LowerBoundInt}
Let $G$ be a finite abelian $p$-group
and let $A_1, \ldots, A_t$ be cyclic subgroups of maximal order and index $s$ in $G$. Let $\omega$ be the number of distinct subgroups $\Omega(A_i)$ where $1\leq i\leq t$.
Then 
\begin{equation}
|\mathcal{X}_G(A_1, \ldots ,  A_t)| \geq s\left(1 - \frac{\omega-1}{p}\right).
\end{equation}
\end{theorem}

\begin{proofof}
We induce on $\omega$. 
If $\omega = 1$ then all $A_i$ share the same subgroup of order $p$. 
Hence $A_1 \cap \ldots \cap A_t > 1$ which implies that 
$|\mathcal{X}_G(A_1, \ldots,  A_t)| = |\mathcal{X}_G(A_1)| = s$ 
and thus our induction begins.

Assume now that $\omega > 1$ and that the claim holds for smaller values of $\omega$. 
Without loss of generality let 
$\Omega(A_u) = \Omega(A_{u+1}) = \ldots = \Omega(A_t)$ and 
$\Omega(A_i)\ne\Omega(A_u)$ for $i < u$. 
The inductive hypothesis implies 
\[
|\mathcal{X}_G(A_1, \ldots , A_{u-1})| \geq s\left(1 - \frac{\omega-2}{p}\right).
\]
From \autoref{Prop:NumberCommonComplementAB} we know 
\begin{align*}
    |\mathcal{X}_G(A_1, \ldots , A_{u-1}) \cup \mathcal{X}_G(A_u, \ldots, A_t)| &\leq |\mathcal{X}_G(A_1) \cup \mathcal{X}_G(A_u)| \\
    &= |\mathcal{X}_G(A_1)|+|\mathcal{X}_G(A_u)|-|\mathcal{X}_G(A_1,A_u)|\\
    &=2s- s\Bigl(1-\frac{1}{p}\Bigr) = s\Bigl(1+\frac{1}{p}\Bigr).
\end{align*}
By inclusion-exclusion, we conclude
\begin{align}
|\mathcal{X}_G(A_1, \ldots , A_{t})| 
&=    |\mathcal{X}_G(A_1, \ldots , A_{u-1}) \cap \mathcal{X}_G(A_u, \ldots, A_t)|\\
&=    |\mathcal{X}_G(A_1, \ldots ,  A_{u-1})| +|\mathcal{X}_G(A_u)| -  
      |\mathcal{X}_G(A_1, \ldots  ,A_{u-1}) \cup \mathcal{X}_G(A_u)| \\
&\geq s\left(1 - \frac{\omega-2}{p}\right) +s -
s\left(1+\frac{1}{p}\right)\\
%&\geq s\left(2 - \frac{\omega-2}{p}\right) -s\left(1+\frac{1}{p}\right)\\ 
&=    s\left(1 - \frac{\omega-1}{p}\right)                    
\end{align}
and the theorem follows.   
\end{proofof}

As an immediate corollary we have the following.
\begin{corollary}\label{Cor:LowerBoundInt}
Let $G$ be a finite abelian $p$-group
and let $A_1, \ldots, A_t$ be cyclic subgroups of maximal order of $G$ then 
$|\mathcal{X}_G(A_1, \ldots  A_t)| \geq s\left(1 - \frac{t-1}{p}\right)$.
\end{corollary}

\section{Proof of \autoref{Thm:B}}

We start with the non-existence part of \autoref{Thm:B}.

\begin{lemma}\label{Lem:noccr}
Let $A$, $B$, $C$ be cyclic subgroups of the finite abelian
group $G$ such that $G_2 = A_2 \times B_2 = A_2 \times C_2 = B_2 \times C_2\ne 1$.
Then $\ccr_G(A, B, C) = \emptyset$.
\end{lemma}

\begin{proofof}
By hypothesis, $\Omega(G_2)=\{1,a,b,c\}$ where $a\in A$, $b\in B$ and $c\in C$. 
Every complement $S$ of $A$ must contain an involution, which lies either in $B$ or in $C$. Hence, $S$ cannot be a common complement of $A,B,C$. 

\end{proofof}

The next theorem addresses the key configuration of \autoref{Thm:B}.

\begin{theorem}\label{Thm:ccrExistencefor3evensubs}
Let $A$, $B$, $C$ be cyclic subgroups of order $2^n$ 
of the abelian $2$-group $G = ABC$. 
Assume further that $A \cap B = A \cap C = 1$ 
and that $|A \cap BC| = 2^m$, $|B \cap C| = 2^k$ 
for some non-negative integers $k$, $m$. 
In case $k > 0$, there is a common complement for
the three subgroups, while if $k = 0$ we have the
following cases:
\begin{enumerate}[label={\upshape(\roman*)}]
\item\label{item:ccr3evensubs1} if $m = 0$, then $A$, $B$, $C$ share a common complement,
\item\label{item:ccr3evensubs2} if $m = n$, 
then there is no common transversal, while
\item\label{item:ccr3evensubs3} if $0 < m < n$, there is no common complement for $A$, $B$, $C$, 
but there is a common transversal. 
\end{enumerate}
\end{theorem}

\begin{proofof}
If $k>0$, then $\Omega(B)=\Omega(C)$ and the claim follows from \autoref{cor:pgrps}.
Now let $k=0$, that is $A \cap B = A \cap C = B \cap C = 1$.

\ref{item:ccr3evensubs1} Here $G = A \times B \times C$ and 
it is easy to see that the subgroup $H = \gensub{ac} \times \gensub{bc}$ 
is a common complement for $A$, $B$, $C$ in $G$,
where $A = \gensub{a}$, $B = \gensub{b}$, $C = \gensub{c}$. 

\ref{item:ccr3evensubs2}
This part follows from \autoref{Lem:noccr}.
% applied to a group with a trivial Hall odd-order subgroup.

\ref{item:ccr3evensubs3}
Finally, let $0<m<n$. 
We assume that a common complement $H$ for $A, B$ and $C$ exists 
and we will derive a contradiction. 
In this case, $G \cong C_{2^n} \times C_{2^n} \times C_{2^r}$ 
with $r = n - m > 0$ and any common complement 
is a group isomorphic to $C_{2^n} \times C_{2^r}$. 
Thus we have $G = H \times A = H \times B = H \times C$ and so 
\[
BC = BC \cap G = BC \cap BH = B \times (BC \cap H)\,.
\]
As $BC = B \times C \cong C_{2^n} \times C_{2^n}$, 
we conclude that $BC \cap H$ is a cyclic group of order $2^n$. 
Now note that $BC$ has only three involutions, 
namely $b^{2^{n-1}}, c^{2^{n-1}}$ and $(bc)^{2^{n-1}}$.
Since $H \cap B = H \cap C = 1$,
we see that the only involution in $BC \cap H$ is $(bc)^{2^{n-1}}$.
Similarly, since $A \cap B = A \cap C = 1$,
the unique involution of 
the non-trivial cyclic subgroup $BC \cap A$ 
is $(bc)^{2^{n-1}}$. 
We conclude that 
\[
(bc)^{2^{n-1}} \in A \cap BC \cap H,
\]
contradicting the fact that $A \cap H = 1$.

To conclude the proof, we must show that $\ccr_G(A, B, C)$ is non-empty.
Since $A \cap BC = \gensub{a^{2^r}}$, we may assume that $a^{2^r}= (bc)^{2^r}$.  

We define the map $\sigma:\{0,\ldots,2^n-1\}\to\{0,\ldots,2^n-1\}$, 
as 
\[
\sigma (i) = t \cdot 2^m + \frac{i-t }{2^r}, 
\]
where $0\leq t<2^r$ with $i\equiv t\pmod{2^r}$.
We note first that $0\leq \frac{i-t}{2^r}<2^m$ and
\[
0\leq \sigma(i) \leq (2^r-1)2^m + \frac{i-t}{2^r} < 2^n.
\]
Suppose next that $\sigma(i) = \sigma(j)$ with $j\equiv t' \pmod{2^r}$. 
Computing modulo $2^m$, 
we obtain $i-t \equiv j-t' \pmod{2^n}$ and hence $i-t = j-t'$. 
But then $t2^m = t'2^m$ and $t = t'$ as well as $i = j$. 
We have shown that $\sigma$ is a permutation.
We also have
$i- \sigma(i) = -t \cdot 2^m + \frac{i(2^r- 1)+t }{2^r} $. Furthermore, if $i \equiv v \pmod{2^r}$ then 
\begin{equation}\label{eq:ccr2gr}
(i -\sigma(i))- (v-\sigma(v) ) = \frac{(i-v)(2^r-1)}{2^r}
\end{equation}
Let $X = \left\{ \, b^i c^{\sigma(i)} \mid i \in \{0,\ldots,2^n-1\} \, \right\}$ and 
\[
Y= \bigcup_{j=0}^{2^r-1} a^jc^{j\cdot 2^m} X.
\]
We claim that $Y \in \ccr_G(A, B, C)$.

It suffices to show that all the elements of $Y$ are distinct
$\pmod{A}$, $\pmod{B}$ and $\pmod{C}$, 
as then we would also have that $|Y| = 2^r \cdot 2^n$. 
\begin{itemize}
\item Let $a^j c^{j 2^m} b^i c^{\sigma(i)} $ and 
$a^u c^{u 2^m} b^v c^{\sigma(v)} \in Y$ 
for some $j, u \in \{0, 1, \ldots, 2^r-1\}$ and $i, v \in [2^n]$. 
Suppose first that 
$a^j c^{j 2^m} b^i c^{\sigma(i)} \equiv a^u c^{u 2^m} b^v c^{\sigma(v)} \pmod{B}$. 
Then 
\[
a^{j-u} \equiv c^{2^m(u-j) +\sigma(v)- \sigma(i)} \pmod{B}.
\]
Hence $a^{j-u} \in A \cap BC$ and thus $2^r \mid j-u$ 
which in turn implies that $j = u$. 
Hence $c^{\sigma(v)- \sigma(i)}\in B$ 
and so $\sigma(v)\equiv \sigma(i) \pmod{2^n}$. 
But $\sigma \in S_{2^n}$, 
so $i$ is necessarily equal to $v$. 
Thus the elements in $Y$ are distinct $\pmod{B}$.

\item Assume now that 
$a^j c^{j 2^m} b^i c^{\sigma(i)} \equiv a^u c^{u 2^m} b^v c^{\sigma(v)} \pmod{C}$ 
which implies that $a^j b^i \equiv a^u b^v \pmod{C}$. 
Hence $a^{j-u} \equiv b^{v-i} \pmod{C}$, 
which in turn yields that $2^r \mid j-u $ and thus $j = u$. 
Therefore $b^{v-i}\in C$ and so $2^n \mid v-i $. 
Thus $v = i$ and the elements of $Y$ are distinct $\pmod{C}$.

\item Lastly, assume that 
$a^j c^{j 2^m} b^i c^{\sigma(i)} \equiv a^u c^{u 2^m} b^v c^{\sigma(v)} \pmod{A}$.
Then $ c^{j 2^m} b^i c^{\sigma(i)} \equiv c^{u 2^m} b^v c^{\sigma(v)} \pmod{A}$ 
and thus $c^{2^m(j-u) +\sigma(i) -\sigma(v) } \equiv b^{v-i} \pmod{A}$. 
So there exists $t \in [2^m]$ such that 
\[
i-v \equiv t 2^r \pmod{2^n} \quad\quad \mbox{and} \quad\quad 
2^m(j-u) +\sigma(i) -\sigma(v) \equiv t 2^r \pmod{2^n}.
\]
We conclude that
\[
2^m(j-u) \equiv (i-\sigma(i)) - (v -\sigma(v)) \pmod{2^n}.
\]
In view of Equation~\eqref{eq:ccr2gr} we get 
\[
2^m(j-u) \equiv \frac{(i-v)(2^r-1)}{2^r} \pmod{2^n}.
\]
But $2^n = 2^{m+r}$ and thus $i-v \equiv 0 \pmod{2^n}$. 
So $i = v$. 
Hence $2^m (j-u) \equiv 0 \pmod{2^n}$ 
and so $2^r \mid j -u$ which yields that $j = u$.
\end{itemize}
The proof is complete.
\end{proofof}

We can now prove that a similar result to that of \autoref{Prop:ReduceToSylow_complements} 
holds for common transversals of three cyclic subgroups.

\begin{corollary}\label{Cor:ReduceToSylow_ccr}
Let $A$, $B$, $C$ be cyclic subgroups of $G = ABC$ with $A \cap B \cap C=1$. 
Then $A$, $B$, $C$ share a common transversal in $G$
if and only if $A_p$, $B_p$, $C_p$ share a common transversal in $G_p$
for all prime divisors $p$ of $|G|$.
\end{corollary}

\begin{proofof}
Assume first that $A_p$, $B_p$, $C_p$ share a common transversal $T_p$ in $G_p$
for all prime divisors $p$ of $|G|$. Then the product $T = \prod_{p} T_p$ is a common transversal of $A$, $B$, $C$ in $G$, as we can easily verify. 
For the other direction, we first note that according to \autoref{Thm:A}, common transversals always exist for three cyclic subgroups of odd order in their product group. Hence we assume that $G_2 \neq 1$ and it suffices to show that if $ \ccr_{G_2}(A_2, B_2, C_2)= \emptyset $ then $\ccr_{G}(A, B, C) =\emptyset $. 
In view of \autoref{Thm:ccrExistencefor3evensubs} we have $\ccr_{G_2}(A_2, B_2, C_2) =\emptyset$ if and only if $A_2 \leq B_2 C_2$ while $B_2 \cap C_2 = 1$ (after some rearrangement of $A$, $B$, $C$). But $A_2, B_2, C_2$ are all cyclic groups of the same order with trivial intersection while their product is $G_2 \neq 1$ and thus $G_2 = A_2 \times B_2= B_2 \times C_2 = A_2 \times C_2$. We are therefore in the situation described in \autoref{Lem:noccr} and thus we get 
$\ccr_{G}(A, B, C) =\emptyset$. 
\end{proofof}

\autoref{Thm:B} is now an easy consequence 
of the preceding corollary.

\begin{proofofB} 
By \autoref{Prop:GisProd}, we may assume that $G=ABC$. Then $A$, $B$ and $C$ are cyclic of maximal order, so they are complemented in $A$. If $|A|$ is odd, then the claim follows from \autoref{Thm:A}.
%If $A$ is of odd order then we clearly have a common transversal 
%for $A$, $B$, $C$ by applying \autoref{Thm:A} 
%along with~\autoref{Prop:GisProd} 
%and~\autoref{Prop:IntersectionIsTriv}. 
We may assume therefore, that $|A|$ is even.
In view of \autoref{Cor:ReduceToSylow_ccr}, 
we have that $\ccr_G(A, B, C) = \emptyset$ 
if and only if $\ccr_{G_2}(A_2, B_2, C_2) = \emptyset$. 
But $\ccr_{G_2}(A_2, B_2, C_2)$ is the empty set 
if and only if $\ccr_{X/Y}(A_2/Y, B_2/Y, C_2/Y) = \emptyset$, 
where $X= A_2 B_2 C_2$ and $Y = A_2 \cap B_2 \cap C_2$. 
Appealing to \autoref{Thm:ccrExistencefor3evensubs} completes the proof.
\end{proofofB}

\section{Some more constructions }\label{Sec:Constructions}
In this section we will provide some methods 
to construct a common transversal for three subgroups $A$, $B$, $C$ of $G$. 
We start with the following generalization of Theorem 1 in~\cite{kolountzakis1997multi}.

\begin{theorem}\label{Thm:AxB_i}
Let $B_1, \ldots, B_t$ be subgroups of $G$ of the same order $m$. 
Let $X = \prod_{i=1}^tB_i $, and assume that $\ccr_X(B_1, \ldots ,B_t) \neq \emptyset$. 
If $A \leq G$ with $|A|= m$ and $A X = A \times X$, 
then $\ccr_{G}(A, B_1, \ldots , B_t) \neq \emptyset$.
\end{theorem}

\begin{proofof}
Clearly, in view of \autoref{Prop:GisProd}, 
we may assume that $G = AX = A \times X$. 
Let $T \in \ccr_X(B_1, \ldots, B_t)$.
Then $|T|= |X|/m$, while $|G|= m \cdot |X|= m^2 \cdot |T|$. 
Observe that for every $b \in B_1$, 
the set $bT \in \ccr_{X}(B_1, \ldots , B_t)$. 
In addition, 
\begin{equation}\label{eq:TrivialIntersOfbT}
bT \cap b'T = \emptyset, \ \text{for all distinct } b, b' \in B_1, 
\end{equation}
or else we would get $bt_1 = b't_2$ for distinct $t_1, t_2 \in T$, 
contradicting the fact that $T$ is a transversal for $B_1$. 
We write $B_1 = \{b_1, \ldots, b_m\}$ and $A = \{a_1, \ldots, a_m\}$ 
and we claim that the set 
\[
D = \bigcup_{i=1}^m a_ib_i T
\]
is an element of $\ccr_{G}(A, B_1, \ldots, B_t)$. 
 
\noindent
We first show that no two elements in $D$ are in the same $A$ or $B_i$-coset 
for all $i \in [t]$. 
Assume first that 
\[
a_ib_it_1 \equiv a_j b_j t_2 \pmod{A}
\]
for $ a_i, a_j \in A, \, \, b_i, b_j \in B_1$ and $ t_1, t_2 \in T$.
Then $b_it_1 \equiv b_j t_2 \pmod{A}$
and so $b_ib_j^{-1} t_1 t_2^{-1} \in A \cap X$. 
As the latter group is trivial, 
we get $b_it_1 = b_j t_2$, $i = j$ and $t_1 = t_2$ by~\eqref{eq:TrivialIntersOfbT}. 

\noindent
Regarding the cosets of $B_i$ for $i \in [t]$ we see that if 
\[
a_ib_it_1 \equiv a_j b_j t_2 \pmod{B_i},
\]
then $a_i a_j^{-1} = b_i^{-1} t_1^{-1} b_j t_2 \in X \cap A$. 
So $a_i = a_j$, that is $i = j$. 
Hence the last congruence implies that $t_1 \equiv t_2 \pmod{B_i}$, 
which means that $t_1 = t_2$, since $T \in \ccr_X(B_1, \ldots, B_t)$.

\noindent
Observe that, as no two elements in $D$ are in the same $A$-coset,
they are necessarily distinct and thus $|D| = m|T| = |X|$ equals the index of $A$ in $G$, 
as well as that of $B_i$ in $G$, for all $i \in [t]$.
We conclude that $D$ is a common transversal for $A, B_1, \ldots, B_t$ in $G$, 
and the theorem follows.
\end{proofof}

As we know, if $B_1$, $B_2$ are subgroups of $G$ of the same order,
then a common transversal exists.
Hence \autoref{Thm:AxB_i} clearly implies the following.

\begin{corollary}
Let $G$ be an abelian group and $A, B, C \leq G$ of the same order $m$. 
Assume further that $ABC = A \times BC$ and let $ T \in \ccr_{BC}(B, C)$.
Then $D = \bigcup_{i=1}^m a_ib_iT $ is a common transversal of $A$, $B$, $C$ in $ABC$, 
where $A = \{a_i\}_{i = 1}^m$ and $B = \{b_i\}_{i = 1}^m$.
\end{corollary}

The above corollary works for any three subgroups $A$, $B$, $C$ of the same order, 
without assuming that they are cyclic, but with the extra hypothesis that $ABC = A \times BC$. 
We have not managed to relax this last hypothesis 
without some restrictions on the type of $A$, $B$, $C$. 
For $A$, $B$, $C$ cyclic subgroups we are able, in some cases, 
to construct the desired common transversal as the next theorem shows.

\begin{theorem}
Let $G = ABC$ with cyclic subgroups $A$, $B$, $C$ of order $p^n$ where $p$ is an odd prime.
Assume further that $A\cap B = A\cap C = B \cap C = 1$ 
while $A \cap BC $ has order $p^{m}$. 
Then for  $a, b, c $ generators of $A$, $B$, $C$ respectively, 
with $a^{p^{n-m}} = (bc)^{p^{n-m}}$  the set 
\[
T = \left\{a^ib^{j-i}c^{-j} : i \in [p^{n-m}], j \in [p^n] \right\} 
\]
is in $\ccr_G(A, B, C)$. 
\end{theorem}

\begin{proofof}
Let $A = \langle a\rangle$ and $r =n-m$. Since $A \cap BC = \gensub{a^{p^r}}$,   we may assume that $a^{p^r}= (bc)^{p^r}$ for appropriate generators $b, c$ of $B$ and $C$. 
We argue that the set $ T = \left\{a^ib^{j-i}c^{-j} : i \in [p^{r}], j \in [p^n]\right\}$ is in $\ccr_G(A, B, C)$.

 \noindent
 Let $a^i b^{j-i} c^{-j}, a^u b^{v-u} c^{-v} \in T$
 for some $i, u \in [p^{r}]$ and $j, v \in [p^n]$.
 Suppose that $a^i b^{j-i} c^{-j} \equiv a^u b^{v-u} c^{-v} \pmod{B}$.
 Then $a^{i-u} \equiv c^{j-v} \pmod{B}$ 
 and thus $a^{i-u} \in A \cap BC$.
 Thus  $ p^{r}$
 divides $i-u$ and so $i = u$  as they are in $[p^r]$.
 So $c^{j-v} \in B \cap C = 1$, 
 which yields that $p^{n}$ divides $j - v$ and thus $j = v$.
 We see therefore that each element in $T$ defines a unique coset of $B$ in $G$.
 The proof that each element in $T$ defines a unique coset of $C$ in $G$
 is entirely analogous so we omit it and we deal next with the case of $A$.
 The congruence here is $b^{j-i} c^{-j} \equiv b^{v-u} c^{-v} \pmod{A}$
 and it yields $b^{(j-v)-(i-u)} \equiv c^{j-v} \pmod{A}$.
 Thus $b^{(j-v)-(i-u)} c^{-(j-v)} \in A \cap BC = \gensub{(bc)^{p^r}}$.
 Since $B \cap C = 1$, it follows that there exists a $t$ such that
 \[
 (j-v) - (i-u) \equiv t p^r \pmod{p^n} \quad\quad \mbox{and} \quad\quad
 v-j \equiv t p^r \pmod{p^n}.
 \]
 Thus $u-i \equiv 2 t p^r \pmod{p^n}$ and so $u-i \equiv 0 \pmod{p^r}$. 
It  follows that $i = u$ as $i, u \in [p^r]$.
 Therefore  $j-v \equiv t p^r \equiv v-j  \pmod{p^n}$ and so $2(j-v) \equiv 0 \pmod{p^n}$. As $p$ is an odd prime we get $p^n \mid j-v$ forcing $j = v$,
 as wanted.

\noindent
We have therefore shown that all the elements of $T$ are in distinct $A, B$ and $C$ cosets. 
Hence the elements of $T$ are pairwise distinct while the cardinality 
$\left\lvert  T \right\rvert  = p^{2n-m}$ 
is the correct one and thus $T \in \ccr_G(A, B, C)$.
\end{proofof}

\noindent 
{\bf  Data availability statement}  Data sharing was not applicable to this article as no datasets were generated or
analyzed.

{\bf  Funding and/or Conflicts of interests/Competing interests.}
The authors have no conflicts of interest to declare that are relevant to the content of this article.

\end{document}